\documentclass[12pt]{amsart}
\usepackage{amscd}
\usepackage{amssymb}

\numberwithin{equation}{section}

\newcommand{\lab}[1]{\label{#1}}

\newtheorem{Thm}{Theorem}[section]
\newtheorem{Prop}[Thm]{Proposition}
\newtheorem{Lem}[Thm]{Lemma}
\newtheorem{Cor}[Thm]{Corollary}

\theoremstyle{remark}
\newtheorem{Rem}[Thm]{Remark}
\newtheorem{Exa}[Thm]{Example}

\theoremstyle{definition}

\newtheorem{Def}[Thm]{Definition}

\newcommand{\bbR}{{\mathbb{R}}}
\newcommand{\bbQ}{{\mathbb{Q}}}
\newcommand{\bbZ}{{\mathbb{Z}}}
\newcommand{\cs}{{\rm CS}}
\newcommand\braket[2]{\left\langle{ {#1}\,,\,{#2} }\right\rangle}
\newcommand{\de}{\partial}

\newcommand{\heta}{\hat\eta}

\newcommand{\GE}{\mathcal{G}_E}
\newcommand{\GEo}{\mathcal{G}_E^0}
\newcommand{\GEi}[2]{\mathcal{G}_{E;{#1},{#2}}}

\newcommand{\dpf}{d_{\pi^{-1}F}}
\newcommand{\wU}{\widetilde{U}}
\newcommand{\weta}{\widetilde{\eta}}

\newcommand{\ord}{\operatorname{ord}}

\newcommand{\Tr}{\operatorname{Tr}}

\newcommand\qq{\rm}

\newcommand\jdg[1]{{\qq J.\ Diff.\ Geom.\ \bf #1}}

\begin{document}
\title{Integral Invariants of 3-Manifolds. II}
\author[R.~Bott]{Raoul Bott}
\address{Department of Mathematics\\
Harvard University\\
Cambridge, MA~02138, USA}
\email{bott@math.harvard.edu}

\author[A.~S.~Cattaneo]{Alberto S.~Cattaneo}
\address{
Lyman Laboratory of Physics\\
Harvard University\\
Cambridge, MA~02138, USA}
\email{cattaneo@math.harvard.edu}

\date{August 28, 1997}

\thanks{We thankfully acknowledge support from NSF for R.B. and
from INFN (grant No.\ 5565/95) 
and from DOE (grant No.\ DE-FG02-94ER25228,
Amendment No.\ A003) for A.S.C}

\begin{abstract}
This note is a sequel to our earlier paper of the same
title \cite{BC} and describes invariants of rational
homology 3-spheres associated to acyclic orthogonal local
systems. Our work is in the spirit of the Axelrod--Singer
papers \cite{AS}, generalizes some of their results,
and furnishes a new setting for the purely topological
implications of their work.
\end{abstract}
\maketitle

\section{Introduction}\lab{sec-intro}
This note is an addendum to our earlier paper of the same title \cite{BC}.
Our aim here will be to construct invariants for framed
3-dimensional homology spheres $(M,f)$, associated to an acyclic
orthogonal local system $E$ on $M$.

Like in our earlier note, we follow the guidelines of the
Axelrod--Singer paper \cite{AS} on the asymptotic of the
Chern--Simons theory, and we have again put aside the
physics inspired aspects of the subject, concentrating
our efforts on the construction of potential configuration-space
integral invariants of $(M,f)$. More precisely we are seeking
invariants that depend on the diffeomorphism type of $M$
and the {\em homotopy class} of the framing $f$.

For simplicity we assume throughout that $M$ is a connected, oriented
3-dim\-en\-sion\-al homology sphere so that---up to
conjugacy---local systems over $M$ are classified by representations
of $\pi_1(M;p)$ where $p$ is some fixed point in $M$.

Our invariants are associated to local systems $E$ which are
induced by an {\em orthogonal} representation $\rho_E$ of $\pi_1(M;p)$
on some $\bbR^m$, and we call such systems orthogonal.
Furthermore, a local system $E$ is called acyclic if
$H^*(M;E)=0$.

With this understood, our principal observation is given by the following
\begin{Thm}
An orthogonal and acyclic local system $E$ over $M$ gives rise
to a purely combinatorial graph cohomology $\GE$, and if\/
$\Gamma\in\GE$ is a connected trivalent cocycle of $\GE$,
then\/ $\Gamma$ determines a numerical invariant $I_\Gamma(M,f)$
of $(M,f)$.
This invariant has the structural form:
\begin{equation}
I_\Gamma(M,f) = A_\Gamma(M) + \phi(\Gamma)\,\cs(M,f),
\end{equation}
where $A_\Gamma(M)$ denotes a sum of configuration-space integrals
specified by $\Gamma$ and a fixed---but arbitrary---Riemannian
structure $g$ on $M$, $\phi(\Gamma)$ is a number universally
associated to\/ $\Gamma$, and\/ $\cs(M,f)$ stands for the Chern--Simons
integral of $M$ relative to $f$ and the Levi-Civita connection of $g$.
\lab{thm-1}
\end{Thm}

Combined with the invariants described in \cite{BC}, where we treated
the trivial coefficient system $\bbR$---which is not
acyclic---, one is therefore in the possession of a large number
of integral invariants of $(M,f)$, and it would be very interesting
to understand their relation to the finite type invariants described
in \cite{FTI} (see also \cite{BGRT})
and whether in their totality they are in any sense exhaustive.

The proof of Theorem \ref{thm-1}, as well as the precise definition
of $\GE$, will be brought in sections \ref{sec-Theta} and 
\ref{sec-hi}, and runs pretty well along the lines of our earlier paper
\cite{BC}. In fact the acyclicity of $E$ allows for a simplification
of the initial step in our procedure, and we will explain this
phenomenon here and now.

Recall that the compactified configuration space $C_2(M)$ is a manifold
with boundary isomorphic to $M\times M$ with its diagonal,
$\Delta$, blown up. Thus $\de C_2(M) = S$ is isomorphic
to the unit sphere bundle of the tangent bundle of $M$.

This situation therefore gives rise to the diagram below
consisting of sections of the exact sequences associated
to the pair $(C_2(M),S)$ and its image $(M\times M,\Delta)$
under the natural projection $\pi$ of $C_2(M)$ to $M\times M$.
{\tiny
\begin{equation}
\begin{CD}
@>>> H^2(C_2(M)) @>>> H^2(S) @>{\delta}>> H^3(C_2(M),S)
@>>> H^3(C_2(M)) @>>> \\
@.  @A{\pi^*}AA @A{\pi^*}AA @A{\pi^*}A{\thickapprox}A @A{\pi^*}AA \\
@>>> H^2(M\times M) @>>> H^2(\Delta) @>{\delta}>> H^3(M\times M,\Delta)
@>>> H^3(M\times M) @>>>
\end{CD} 
\lab{seq}
\end{equation} }

The vertical isomorphism in the third column of \eqref{seq}
follows from excision near the blow-up of $\Delta$.
Note also that this diagram is acted upon by the involution $T$
which exchanges the factors in $M\times M$, and each of the sequences
therefore splits canonically into a $+$ and $-$ part
corresponding to the $\pm1$ eigenvalues of $T$.
In the bottom row $H^*(\Delta)$ is clearly invariant under $T$
so that the antisymmetric part of \eqref{seq} reduces to
\begin{equation}
\begin{CD}
 H^2_-(C_2(M)) @>>> H^2_-(S) @>{\delta}>> H^3_-(C_2(M),S)\\
  @. @. @A{\pi^*}A{\thickapprox}A \\
  @. 0 @>>> H^3_-(M\times M,\Delta)
@>{\thickapprox}>> H^3_-(M\times M)
\end{CD}
\lab{redseq}
\end{equation} 

Now in \cite{BC} we showed that the form $\eta$ given by half
the Euler form of the tangent bundle along the fiber in the
fibering $S\to\Delta$ generates $H^2(S;\bbR)$
as a module over $H^*(\Delta;\bbR)$ and that this $\eta$
is antisymmetric: $T^*\eta=-\eta$. In short $[\eta]$
generates $H^2_-(S;\bbR)$.

So far our discussion involved the constant coefficient system $\bbR$.
But the sequences \eqref{seq} and \eqref{redseq} as well as the
action of $T$ remain valid for a general local system $E$ on $M$,
provided we use the local system
$F = \pi_1^{-1}E\otimes\pi_2^{-1}E$
on $M\times M$, and $\pi^{-1}F$ on $C_2(M)$.

This understood, assume now that $E$ is orthogonal and acyclic.

The orthogonality gives rise to an arrow
\begin{equation}
I:\bbR\to E\otimes E,
\end{equation}
defined by sending 1 to $\sum_i e_i\otimes e_i$, where
$\{e_i\}$ is any orthonormal frame in $E$. We may therefore
also apply $I$ to $\eta$ to obtain a {\em closed}
form $I(\eta)\in\Omega^2(S;E\otimes E)$.

Next observe that the K\"unneth formula implies that
\[
H^*(M\times M;\pi_1^{-1}E\otimes\pi_2^{-1}E)=
H^*(M;E)\otimes H^*(M;E).
\]

Hence, under our acyclicity assumption, all the terms on the right
of $\delta$ in \eqref{redseq} vanish! It follows immediately
that the class of $I(\eta)\in\Omega^2(S;E\otimes E)$
is in the image of a class $[\heta]\in H^2(C_2(M);\pi^{-1}F)$.

Actually we need the following slight refinement of this assertion:
\begin{Lem}
Under the orthogonality and acyclicity assumptions
there exists a form $\heta\in\Omega^2(C_2(M);\pi^{-1}F)$
with the following properties:
\begin{enumerate}
\item The restriction of $\heta$ to $S$ is $I(\eta)$:
$i_\de^*\heta=I(\eta)$.
\item $\heta$ is closed under $\dpf$: $\dpf\heta=0$.
\item $\heta$ is antisymmetric: $T^*\heta=-\heta$.
\end{enumerate}
\lab{lem-heta}
\end{Lem}
The construction of $\heta$ proceeds precisely along the guidelines
given in \cite{BC}.

Let $U$ be a tubular neighborhood of $\Delta$ in $M\times M$,
and let $p:U\to\Delta$ be a projection which
fibers $U$ over $\Delta$
into discs on which $T$ acts linearly as the antipodal map,
and such that $\de U$ can be identified with $S$.
Then $\wU=\pi^{-1}U$ has the structure of $S\times[0,1]$
and hence fibers over $S$ with the unit interval as fiber.
We write $\sigma:\wU\to S$ for the projection onto $S$
in this fibering, and note that $T$ acts on $S\times[0,1]$
by the antipodal map on $S$ crossed with the identity on $[0,1]$.
Now choose a smooth function $\chi$ on $[0,1]$ which is
identically $+1$ near 0 and identically 0 near $+1$, and write
$\chi$ also for its pullback to $\wU$.
It follows that the form
\[
\weta = \sigma^*I(\eta)\,\chi
\]
on $\wU$ extends by $0$ to a form on all of $C_2(M)$ with values
in $\pi^{-1}F$. It is also clear that $\weta$ restricts to $I(\eta)$
on $S$, so that $\dpf\weta$ represents $\delta(I(\eta))$
in the upper sequence. On the other hand $\dpf\weta$ vanishes
identically near $S$ and so may be considered an antisymmetric
form on $M\times M$. But then, by the acyclicity assumption,
there must exist an antisymmetric form $\alpha$ on $M\times M$
such that
\[
d_F\alpha = \dpf\weta.
\]
Now $\heta=\weta-\pi^*\alpha$ has all the desired properties.

\section{The $\Theta$-invariant} \lab{sec-Theta}
Using the closed form $\heta$ defined in the previous section,
see Lemma \ref{lem-heta}, we can define an invariant
for the framed 3-dimensional homology sphere $(M,f)$.

First we notice that $\heta^3$ is a 6-form on the 6-dimensional
space $C_2(M)$ which takes values in 
$\pi_1^{-1}E^{\otimes3}\otimes\pi_2^{-1}E^{\otimes3}$.
If we  associate
to each vertex $i$ ($i=1,2$) a homomorphism
\begin{equation}
\rho_i:\bbR\to E\otimes E\otimes E
\end{equation}
which is equivariant as a module over $\pi_1(M)$,
then we obtain the {\em closed} real-valued
6-form $\braket{\rho_1\rho_2}{\heta^3}$.
Here $\braket{\cdot}{\cdot}$ denotes the scalar product 
over $E$ and its extensions to $E^{\otimes3}$ and to 
$\pi_1^{-1}E^{\otimes3}\otimes\pi_2^{-1}E^{\otimes3}$.

\begin{Rem}
The existence of such homomorphisms depends on $E$. In some
cases, the only possible choice will be the trivial one: $\rho=0$. 

If the vector
space spanned by these homomorphisms has dimension greater than one,
then one can choose $\rho_1$ and $\rho_2$ linearly independent.
\end{Rem}
\begin{Exa}
A particular case, considered in \cite{AS}, occurs when
$E$ is the adjoint representation of a compact Lie group $G$. 
Then a natural
choice for $\rho$ is obtained by using the structure constants
$f_{abc}$ relative to a left- and right-invariant inner product
on the Lie algebra of $G$; namely,
\[
\rho(x) = x\,\sum_{abc} f_{abc}\;e_a\otimes e_b\otimes e_c.
\]
The equivariance under the full group $G$ ensures the equivariance
under the action of $\pi_1(M)$. Notice that the antisymmetry of
the structure constants implies that this homomorphism is completely
antisymmetric.

Note that if $E$ denotes a representation of such a Lie group $G$,
the equivariant homomorphisms are dual
to the projections to the trivial representations
in $E\otimes E\otimes E$. Again,
the equivariance under the full group $G$ ensures the equivariance
under the action of $\pi_1(M)$.
\lab{exa-adj}
\end{Exa}

\begin{Exa}
If $E_j$ denotes the irreducible representation of spin $j\in\bbZ/2$ of
$SU(2)$, then the Clebsch--Gordan formula,
\[
E_j\otimes E_k = \bigoplus_{l=|j-k|}^{j+k} E_l,
\]
implies that
\[
E_j\otimes E_j\otimes E_j =
\bigoplus_{k=j-\lfloor j\rfloor}^j(2k+1)E_k \oplus
\bigoplus_{k=j+1}^{3j}(3j-k+1)E_k.
\]
So $E_j^{\otimes3}$ contains no trivial representations
if $j$ is a half-integer and one trivial representation if
$j$ is an integer. In the case $j=1$ we recover the choice of
example \ref{exa-adj}. Notice that this projection is
obtained by selecting the representation of spin $j$ in
$E_j\otimes E_j$, then by tensoring by the last copy of $E_j$, 
and finally by projecting on the trivial representation. Therefore,
all these projections (and the corresponding homomorphisms) are
completely antisymmetric.
\end{Exa}

\begin{Exa}
With the notations of the previous example, consider
$E=E_{1/2}\oplus E_1$. It turns out that $E^{\otimes3}$
contains three copies of the trivial representation:
the first is obtained by choosing the trivial representation
in $E_1^{\otimes3}$; the second by choosing
the trivial representation in 
$E_{1/2}\otimes E_{1/2}\otimes E_1$, and the other two by
cyclic rotations of the second. Notice that the second
projection is obtained by selecting the representation of
spin 1 in $E_{1/2}\otimes E_{1/2}$, then by tensoring by
$E_1$, and finally by projecting on the trivial representation.
Therefore, this projection is symmetric with respect
to the exchange of the spin-$(1/2)$ components.
\lab{exa-121}
\end{Exa}

Integrating the closed form we have obtained
over $C_2(M)$ yields the number
\begin{equation}
A_{(\Theta,\rho_1,\rho_2)} \doteq \int_{C_2(M)} 
\braket{\rho_1\rho_2}{\heta^3},
\end{equation}
which is our first potential invariant. We recall, see \cite{BC},
that the definition of $\eta$ relies on the choice of a metric
on $M$ and of a compatible connection; 
moreover, the construction of $\heta$ requires the choice of a function
$\chi$ and of a 2-form $\alpha$
as explained after Lemma \ref{lem-heta}. An invariant must
be independent of all these choices. Actually, we have
the following
\begin{Thm}
Given a section $f$ of the orthonormal frame bundle, the
combination
\begin{equation}
I_{(\Theta,\rho_1,\rho_2)}(M,f) = 
A_{(\Theta,\rho_1,\rho_2)}(M) -
\frac{\braket{\rho_1}{\rho_2}}4\, \cs(M,f),
\lab{ITheta}
\end{equation}
is independent of all the choices involved (except for the framing).
Here
\begin{multline}
\cs(M,f)=-\frac1{8\pi^2}\int_Mf^*\Tr\left({
\theta\,d\theta +\frac23\,\theta^3}\right)=\\
\frac1{4\pi^2}\int_M f^*\left({
\theta^id\theta_i-\frac13\,\epsilon_{ijk}\theta^i\theta^j\theta^k
}\right),
\end{multline}
is the Chern--Simons integral of the same metric
connection used to define $\eta$.

Thus, $I_{(\Theta,\rho_1,\rho_2)}(M,f)$ 
is an invariant for the framed rational homology
sphere $(M,f)$.
\lab{thm-Theta}
\end{Thm}

\begin{Rem}
In the case discussed in example \ref{exa-adj}, we have
\[
\braket{\rho_1}{\rho_2} = \sum_{abc} f_{abc}f_{abc} =
-c_v\,\dim G,
\]
where $c_v$ is the Casimir of the adjoint representation of $G$.
\end{Rem}

\begin{proof}
As in \cite{BC}, we introduce the unit interval $I$ as
a parameter space, and recall that, as shown there,
letting $\theta$ vary on $I$ corresponds to defining 
on $S\times I$ a form---which we still denote by 
$\eta$---given by half the Euler form of the tangent bundle
along the fiber in the fibering $S\times I\to\Delta\times I$. 

Then all the arguments contained in section \ref{sec-intro} are
still true if we multiply by $I$ each space involved (say, $M$,
$M\times M$, $\Delta$, $C_2(M)$ and $S$), 
since $H^n(I)=\delta_{n0}\,\bbR$. 
In particular, we have a form---which we keep denoting by 
$\heta$---which satisfies the properties of Lemma \ref{lem-heta}
with $C_2(M)$ replaced by $C_2(M)\times I$. (To be precise,
by $\pi$ now we mean the projection 
$C_2(M)\times I\to M\times M\times I$.)

If we denote by $\sigma$ the projection $C_2(M)\times I\to I$,
then
\[
A_{(\Theta,\rho_1,\rho_2),\tau} = 
\sigma_*\braket{\rho_1\rho_2}{\heta^3}
\]
is a function depending on the parameter $\tau\in I$, in whose
variations we are interested. 

To do so, we recall that, given two
spaces $M_1$ and $M_2$ and projections $\pi_i:M_1\times M_2\to M_i$,
Stokes' theorem can be rewritten as
\begin{equation}
d\pi_{i*}\omega = \pi_{i*}d\omega - 
(-1)^{\deg\pi_{i*}^\de\omega}\,\pi_{i*}^\de\omega,
\lab{dpi*}
\end{equation}
where $\pi_{1*}^\de$ denotes integration along the boundary of
$M_2$ and vice versa. (Notice that the signs in \eqref{dpi*}
are correct if integration acts from the right.)

Since $\braket{\rho_1\rho_2}{\heta^3}$ is a closed form, we simply
have
\[
dA_{(\Theta,\rho_1,\rho_2),\tau} = 
\sigma_*^{\de}\braket{\rho_1\rho_2}{\heta^3} =
\braket{\rho_1}{\rho_2}\sigma_*^{\de}\eta^3,
\]
the last identity following from property 1 in Lemma \ref{lem-heta}.

Now we recall that in \cite{BC} (see Lemma 3.15 there)
we proved that
\[
\pi_*^{\de}\eta^3 = \frac14\, p_1,
\]
where $\pi^\de$
is the projection $S\times I\to\Delta\times I$, and
$p_1$ is the first Pontrjagin form on 
$\Delta\times I=M\times I$.

Denoting by $\sigma_M$ the projection
$M\times I\to I$, we finally get
\[
dA_{(\Theta,\rho_1,\rho_2),\tau} =
\frac{\braket{\rho_1}{\rho_2}}4\,\sigma_{M*} p_1,
\]
from which the theorem follows.
\end{proof}

\section{The higher invariants}\lab{sec-hi}
Using the natural projections $\pi_{ij}:C_n(M)\to C_2(M)$
we can pull back the form $\heta$ defined in section \ref{sec-intro}.
We will write
\[
\heta_{ij} = \pi_{ij}^*\heta,
\]
and by property 3 of Lemma \ref{lem-heta} we have
\[
\heta_{ij} = -\heta_{ji}.
\]

These forms on $C_n(M)$ allow for
writing other invariants of the 3-dimensional
homology sphere $M$ associated to cocycles in an appropriate graph 
cohomology (depending on the bundle $E$).

\begin{Def}
We call a decorated graph a graph with oriented and numbered edges and 
numbered vertices (by convention we start the
enumeration by 1). We require edges always to connect distinct
vertices. If two vertices are connected by exactly one edge,
we call that edge {\em regular}.

The edge numbering induces a numbering of the $v_i$ 
half-edges at each vertex $i$,
corresponding to which we attach a homomorphism
\[
\rho_i : \bbR\to E^{\otimes v_i}
\]
which is equivariant as a module over $\pi_1(M)$.

Denoting by $V$ the number of vertices and
by $E$ the number of edges, we grade the collection of 
decorated graphs by
\begin{equation}
\begin{split}
\ord\Gamma &= E-V,\\
\deg\Gamma &= 2E -3V.
\end{split}
\lab{ord}
\end{equation}

\end{Def}

\begin{Rem}
Compared with the decorated graphs we introduced in \cite{BC},
this definition adds two further decorations: the numbering
of the edges and the equivariant homomorphisms attached to the vertices.
\end{Rem}

\begin{Rem}
A trivalent diagram has degree zero, and its order is
given by $m=V/2=E/3$.

We thank S.~Garoufalidis for pointing out that our choice
of the words ``order" and ``degree" is a bit unfortunate, for
people working with finite type invariants
call $m$ the degree (instead of the order)
of a trivalent graph.

However, we prefer to stick to our old notation \cite{BC} since
the term degree is consistent with the cohomology defined
by the coboundary operator $\delta$ (see Proposition 
\ref{prop-delta}).
\end{Rem}

Denoting by $v(\Gamma)$ the set of vertices and by $e(\Gamma)$ the
ordered set of oriented edges in $\Gamma$, we can associate to
the 3-dimensional homology sphere $M$ and to the {\em trivalent}
decorated graph $\Gamma$ the number
\begin{equation}
A_\Gamma(M) \doteq \int_{C_n(M)} 
\braket{\prod_{i\in v(\Gamma)}\rho_i}
{\prod_{(ij)\in e(\Gamma)}\heta_{ij}},
\lab{defAGamma}
\end{equation}
where $n=2\ord\Gamma$ is the number of vertices and $(ij)$ denotes
the edge connecting the vertex $i$ to the vertex $j$.

Next we give the collection of 
decorated graphs the structure of an algebra over 
$\bbQ$ (the product simply being the disjoint union
of graphs). We will denote this algebra by
$\GEo$ and will extend \eqref{defAGamma} by linearity. 

In view of the definition of $A_\Gamma(M)$, we introduce
the following equivalence relation on $\GEo$: 
if two decorated graphs $\Gamma$ and $\Gamma'$ 
differ only by a permutation of order $p$ in
the vertex numbering and by $l$ edge-orientation reversals,
we set
\begin{equation}
\Gamma=(-1)^{(p+l)}\,\Gamma'.
\lab{eqrel}
\end{equation}
Notice that to equivalent
graphs we associate the same number $A_\Gamma(M)$. We will denote
by $\GE$ the algebra of graphs modulo the above equivalence
relation.

Then we introduce an operator $\delta$ on $\GEo$
that acts by contracting a regular edge one at a time in $\Gamma$,
followed by a consistent renumbering of edges and vertices.
To the contraction of the regular edge connecting the vertex $i$ to
the vertex $j$ we associate a sign $\sigma(i,j)$
defined by
\begin{equation}
\sigma(i,j) = \begin{cases}
(-1)^j &\text{if $j>i$,}\\
(-1)^{i+1} &\text{if $j<i$}.
\end{cases}
\lab{sigmaij}
\end{equation}
Assuming that this edge corresponds to the $k$th of the $v_i$ half-edges
at $i$ and to the $l$th of the $v_j$ half-edges at $j$, we attach
to the vertex obtained after contraction the equivariant homomorphism
\[
\widetilde{\rho}_i:\bbR\to E^{\otimes(v_i+v_j-2)}
\]
defined by
\[
\widetilde{\rho}_i = m_{k,v_i+l}(\rho_i\otimes\rho_j),
\]
where $m_{rs}$ denotes the scalar product between the $r$th and
the $s$th terms in the tensor product.

Notice that the homomorphism attached to a vertex after contracting
two different (regular) edges starting from it does not depend on
the order of the contractions. Therefore, the argument we gave
in \cite{BC} is enough to prove the following
\begin{Prop}
The operator $\delta$ descends to $\GE$ and satisfies
$\delta^2=0$ there. Moreover, if we denote by $\GEi nt$
the (equivalence classes of) decorated graphs of order $n$ and
degree $t$, we have
\[
\delta : \GEi nt \to \GEi n{t+1}.
\]
\lab{prop-delta}
\end{Prop}

We call a cocycle an element of the kernel of $\delta$ in
$\GE$. Notice that the action of $\delta$ can be
restricted to the algebra of (equivalence classes of) decorated
connected graphs. Now we are in a position to prove Theorem
\ref{thm-1}.

\begin{proof}[Proof of Theorem \ref{thm-1}]
As in the proof of Theorem \ref{thm-Theta} we introduce
the unit interval $I$ as a parameter space, and define
$\heta$ as a form on $C_n(M)\times I$. Denoting by $\sigma$
the projection $C_n(M)\times I\to C_n(M)$, we define
\[
A_{\Gamma,\tau}(M) = \sigma_*
\braket{\prod_{i\in v(\Gamma)}\rho_i}
{\prod_{(ij)\in e(\Gamma)}\heta_{ij}},
\]
and consider its change as $\tau$ varies on $I$.

Since the integrand form is closed, equation \eqref{dpi*}
implies that $dA_{\Gamma,\tau}(M)$ is given by boundary 
contributions only.

But on the boundary $\heta$ reduces to $I(\eta)$---see Lemma 
\ref{lem-heta}---so we can use essentially the same arguments
as in the trivial coefficient case \cite{BC}. Therefore, we will only
give a brief sketch of the proof here and refer to \cite{BC}
for further details.

First recall that the face
in $\de C_n(M)$ corresponding to the collapse of
$q$ points has the structure of a fibering over $C_{n-q+1}$,
with the $(3q-4)$-dimensional fiber isomorphic to
$C_q(\bbR)$ modulo translations and scalings.

If we denote by $e$ the number of edges connecting two collapsing
vertices, we see that the vertical
form-degree is given by $2e$. Moreover, since we are considering
trivalent graphs, we have the relation $2e+e_0=3q$, where
$e_0$ denotes the number of edges connecting a collapsing
vertex with a non collapsing one.  Therefore, the push-forward
along the fiber of the form associated to the edges connecting 
collapsing vertices yields a form of degree $4-e_0$ if
$e_0\leq4$, and zero otherwise. 

By a theorem due to Axelrod and Singer \cite{AS}---which we
have recast in \cite{BC} in a form suitable to our 
construction---this form must be the pullback of a 
multiple of a characteristic form on $M\times I$, namely, the
constant function or the first Pontrjagin form $p_1$.

The former case corresponds to $e_0=4$, and it can be shown
that the only case when the integral does not vanish is when
this is given by the collapse of just two vertices.
These boundary terms are then taken care of by the requirement
that $\Gamma$ be a cocycle. 

The latter case corresponds to $e_0=0$, that is, to the case
when all points collapse since we assume the diagrams to be connected.
These terms are then taken care of by the correction 
$\phi(\Gamma)\,\cs$.
\end{proof}

\begin{Rem}
It is clear from the above proof that $\phi(\Gamma)$ is linear.
Moreover, if $\Gamma$ is a decorated
graph, we can write
\[
\phi(\Gamma) = \rho(\Gamma)\,\phi_0(\Gamma),
\]
where $\rho(\Gamma)$ is a purely algebraic factor obtained
from the homomorphisms $\rho_i$
by associating a scalar product to each edge in $\Gamma$,
while
$\phi_0(\Gamma)$ is given by the boundary integral involving
the forms $\eta$, and so it is the same as in the trivial coefficient
case.

If $\ord\Gamma$ is even, there exists an orientation reversing
involution under which the integrand form turns out to be odd
(see \cite{AS} or \cite{BC}). Therefore, in this case,
$\phi_0(\Gamma)$ vanishes, and so does $\phi(\Gamma)$.
\end{Rem}

{}From  the above remark and from
Theorems \ref{thm-1} and \ref{thm-Theta},
we get the following
\begin{Cor}
If\/ $\Gamma$ is a connected trivalent cocycle of $\GE$,
and $\rho_1$ and $\rho_2$ are equivariant homomorphisms such that
$\braket{\rho_1}{\rho_2}\not=0$,
then the quantity
\[
J_{\Gamma;\rho_1,\rho_2}(M) 
= A_\Gamma(M) + \frac4{\braket{\rho_1}{\rho_2}}
\,\phi(\Gamma)\, A_{(\Theta,\rho_1,\rho_2)}(M)
\]
is an invariant for the rational homology 3-sphere $M$. 
Moreover, if\/ 
$\ord\Gamma$ is even, then $\phi(\Gamma)=0$. 
\end{Cor}

\section{Discussion}\lab{sec-disc}
The graph cohomology introduced in the previous section
is in principle more general than those introduced in
\cite{AS} and \cite{BC} and might give rise to more general
invariants.

In the case when $E$ is the adjoint bundle of a Lie group
$G$, we can choose all the equivariant homomorphisms 
associated to a trivalent graph to be
determined by the structure constants, as explained in example
\ref{exa-adj}. The cocycles in this case are those studied in
\cite{AS} and come naturally from perturbative Chern--Simons
theory. Antisymmetry of the structure constants implies
that it is enough to give a {\em cyclic order} of the three half-edges
at each vertex. The Jacobi identity then implies that the cocycles
satisfy the so-called IHX relation (see \cite{B-N}).
It is a non-trivial fact (and we thank S.~Garoufalidis for pointing this out)
that these cocycles are in one-to-one correspondence with
the cocycles of the trivial coefficient case.

If $E_{1/2}\oplus E_1$ (where $E_{1/2}$ and $E_1$
denote the representation of $SU(2)$ with spin $1/2$ and 1)
is an orthogonal and acyclic local system, we can choose each
homomorphism as the dual of the second projection
considered in example \ref{exa-121}.
In this case we can think of the diagram as carrying
spin $1/2$ over two of the three half-edges at each vertex
and spin 1 over the last half-edge.
Since each of these
homomorphisms is {\em symmetric} with respect to the
exchange of the two spin-$1/2$ representations, the diagram
is symmetric under the exchange of the corrresponding
half-edges.

It would be interesting to see if this or more general choices of 
the bundle $E$ and of the equivariant homomorphisms give rise
to new inequivalent cocycles.

\section*{Acknowledgements}
We again thankfully acknowledge helpful conversations with
Scott Axelrod, Robin Forman and Cliff Taubes.
We are especially thankful to Stavros Garoufalidis
for pointing out the isomorphisms in various graph
cohomologies mentioned in section \ref{sec-disc}.
One of these isomorphisms was also explained earlier
by Dylan Thurston in his senior thesis at Harvard.


\thebibliography{9}
\bibitem{AS} S. Axelrod and I. M. Singer, ``Chern--Simons Perturbation 
Theory,'' in {\em Proceedings of the XXth DGM Conference}, edited by
S.~Catto and A.~Rocha (World Scientific, Singapore, 1992), 
pp.\ 3--45; ``Chern--Simons Perturbation 
Theory.\ II,'' \jdg{39} (1994), 173--213.
\bibitem{B-N} D.~Bar-Natan, ``On the Vassiliev Knot Invariants,''
Topology {\bf 34} (1995), 423--472.
\bibitem{BGRT} D.~Bar-Natan, S.~Garoufalidis, L.~Rozansky and
D.~P.~Thurston, ``The \AA{}rhus Invariant of Rational Homology
3-Spheres: A Highly Non Trivial Flat Connection on $S^3$,''
q-alg/9706004.
\bibitem{BC} R. Bott and A. S. Cattaneo,
``Integral Invariants of 3-Manifolds,'' dg-ga/9710001, to appear in
J.\ Diff.\ Geom.
\bibitem{FTI} T.~Q.~T.~Le, J.~Murakami and T.~Ohtsuki, ``On a Universal
Quantum Invariant of 3-Manifolds,'' q-alg/9512002,
to appear in Topology.
\end{document}